\input amstex
\documentstyle{amsppt}
\input epsf.tex
\input label.def
\loadbold
\input xypic

\def\C{{\Bbb C}}

\def\CCC{\Bbb O}
\def\R{{\Bbb R}}
\def\Z{{\Bbb Z}}

\def\N{\Cal N}

\def\H{\Cal H}
\def\D{\Delta}
\def\DD{\Bbb D}
\def\DO#1{\D_{\CCC}^{#1}}
\def\G{\Gamma}
\def\T{a}

\def\OO{X}

\def\LLL{\Lambda}

\def\PH{\Phi_{\roman{aut}}}
\def\P#1{\Bbb P^{#1}}
\def\O#1#2{\Bbb O_{#1\, #2}}

\def\a{\alpha}
\def\b{\beta}

\def\e{\varepsilon}

\def\Rp#1{\R\roman P^{#1}}

\def\conj{\operatorname{conj}}

\def\Aut{\mathop{\roman{Aut}}\nolimits}

\def\Sing{\operatorname{Sing}}
\def\Arf{\operatorname{Arf}}

\def\sign{\roman{sign}}

\def\lk{\roman{lk}}

\def\Sym{\operatorname{Sym}}

\def\til{\widetilde}

\let\ge\geqslant
\let\le\leqslant
\def\sm{\smallsetminus}
\let\+\sqcup
\def\oo{\varnothing}

\nologo
\NoBlackBoxes
\leftheadtext{S.~Finashin}
 \topmatter
\title
Deformation classes of real Cayley M-octads
\endtitle
\author Sergey Finashin
\endauthor
\address Middle East Technical University,
Department of Mathematics\endgraf Ankara 06800 Turkey
\endaddress
\abstract\nofrills
We study 8-point configurations in the real projective space forming an
intersection locus of three quadrics and containing no coplanar quadruples.
We found that there exists precisely 8 mirror-pairs of deformation classes of such configurations. 
We describe also the mutual position of these 8 pairs and find the real monodromy
groups acting on the 8-point configurations, for each deformation class.
\endabstract
\endtopmatter

\document

\section{Introduction}
\subsection{Cayley octads}
A {\it Cayley octad} $\OO\subset\P3$ is an 8-point configuration obtained as
the intersection locus of three quadric surfaces $\OO=Q_0\cap Q_1\cap Q_2$.
We allow {\it multiple points} of $\OO$, which appear if the intersection is
not transverse, in which case 8 is the sum of multiplicities.
A simple analysis shows that for any Cayley octad $\OO$
the net of quadrics
$$\N=\{Q_t=t_0Q_0+t_1Q_1+t_2Q_2\},\  t=[t_0:t_1:t_2]\in P^2,$$
is the complete linear system of quadrics passing through $\OO$, and so,
$\N$ and $\OO$ determine each other.

Singular quadrics $Q_t\in\N$ are parameterized by a quartic curve $\H\subset\N\cong P^2$
called the {\it Hessian curve} ({\it Hessian quartic}): it is defined by the determinant
of the symmetric $4\times 4$-matrix defining $Q_t$.
%
The following conditions are known to be equivalent (cf., \cite{Dolgachev, Sect.\,6.3.2}, or \cite{GH, Lemma 6.4}):
\roster\item
The Hessian curve $\H$ associated to a Cayley octad $\OO$ is non-singular.
\item
The Cayley octad $\OO$ contains neither multiple points,
nor coplanar subsets of more than 3 points.
\item
Each quadric $Q_t$ from the net $\N$ associated to $\OO$ is irreducible, and
the points of $\OO$ are non-singular on $Q_t$ for all $t\in P^2$.
\endroster

A Cayley octad $\OO$ is called {\it regular} if these equivalent conditions are satisfied, and
{\it singular} otherwise.

The reality condition (invariance under the complex conjugation)
for Cayley octad $\OO$ and for the net of quadrics $\N$ are obviously equivalent and imply
reality of the Hessian curve $\H$.
We say that a real Cayley octad $\OO$ is {\it maximal} or {\it M-octad} if its eight points are all real.
If $\OO$ contains $k>0$ pairs of conjugate complex points and $8-2k$ real ones, we call it $(M-k)$-octad.

\subsection{Principal Results}
Our principal aim is to enumerate the {\it deformation classes} that is path-connected components
in the space of regular M-octads and to describe the mutual position of these classes.

Our first principal result says that there exist precisely
16 deformation classes of regular M-octads, which can be grouped into 8 pairs
of {\it mirror partner classes} (see Theorems \ref{coarse-classification} and \ref{chirality-of-octads}).
Such pairs of classes, $[\OO]$ and $[\bar\OO]$, are respresented
by a regular M-octad $\OO$ and its image $\bar\OO$ under an orientation reversing projective transformation of $\Rp3$.
The union $[\OO]\cup[\bar\OO]$
is called the {\it coarse deformation class} of an M-octad $\OO$.

The 8 pairs of deformation classes are named $(\O{\a}{\b}^+,\O{\a}{\b}^-)$, where $(\a,\b)\in\{0,2,4\}\times\{0,3,4\}\sm\{(4,3)\}$,
and the course deformation classes are $\O{\a}{\b}=\O{\a}{\b}^+\cup\O{\a}{\b}^-$.
The meaning of indices $\a$ and $\b$ can be explained in terms of decorated graph $\G_\OO$ related to some degenerations of $\OO$.
This graph has the vertex set $\OO$ and its edges are line segments
joining the pairs of vertices, which can be merged by a real variation of $\OO$ formed by regular M-octad.
More precisely, among the two line segments in $\Rp3$ joining a pair of points of $\OO$ we select the one homotopic to
the trace of these points under the merging variation. It follows that an edge with endpoints $x,y\in\OO$ cannot cross
any of the 20 planes passing though the triples of points of $\OO\sm\{x,y\}$.

The combinatorial types of such graphs are presented on Fig.\, 1, Table 1.

\midinsert
\topcaption{Fig.\,1. Graphs $\G_\OO$ and monodromy groups $\Aut_\R(\OO)$ of regular M-octads $\OO\in\O{\a}{\b}$}
\endcaption
\hbox{\hsize78mm\vtop{\hbox{\hskip4mm Table 1. Combinatorial types of $\G_\OO$ for $\OO\in\O{\a}{\b}$}
$\boxed{\xymatrix@1@=12pt{
&\bullet&\bullet\\
\bullet\ar@{-}[ur]^{\b}&&&\bullet\ar@{-}[ul]_{\b}\\
\bullet\ar@{-}[dr]_\b&{}\ar@{}^{\O04}&&\bullet\ar@{-}[dl]^\b\\
&\bullet&\bullet}}
\hskip1pt
\boxed{\xymatrix@1@=12pt{
&\bullet\ar@{-}[r]_\a&\bullet\\
\bullet\ar@{-}[ur]^\b&&&\bullet\ar@{-}[ul]_\b\\
\bullet\ar@{-}[dr]_\b&\ar@{}^{\O24}&&\bullet\ar@{-}[dl]^\b\\
&\bullet\ar@{-}[r]^\a&\bullet}}
\hskip1pt
\boxed{\xymatrix@1@=12pt{
&\bullet\ar@{-}[r]_\a&\bullet\\
\bullet\ar@{-}[ur]^\b\ar@{-}[d]^\a&&&\bullet\ar@{-}[d]_\a\ar@{-}[ul]_\b\\
\bullet\ar@{-}[dr]_\b&\ar@{}^{\O44}&&\bullet\ar@{-}[dl]^\b\\
&\bullet\ar@{-}[r]^\a&\bullet}}$
\vskip0mm
$\boxed{\xymatrix@1@=12pt{
&\bullet\ar@{-}[dl]_\b&\bullet\ar@{-}[dr]^\b\\
\bullet&&&\bullet\\
\bullet\ar@{-}[dr]_\b&\ar@{}^{\O03}&&\bullet\\
&\bullet&\bullet}}
\hskip1pt
\boxed{\xymatrix@1@=12pt{
&\bullet&\bullet\ar@{-}[dr]^\b\\
\bullet\ar@{-}[ur]^\b&&&\bullet\ar@{-}[d]_\a\\
\bullet\ar@{-}[dr]_\b&\ar@{}^{\O23}&&\bullet\\
&\bullet\ar@{-}[r]^\a&\bullet}}$
\vskip0mm
$\boxed{\xymatrix@1@=12pt{
&\bullet&\bullet\\
\bullet&&&\bullet\\
\bullet&{}\ar@{}[r]^{\O00}&&\bullet\\
&\bullet&\bullet}}
\hskip1pt
\boxed{\xymatrix@1@=12pt{
&\bullet\ar@{-}[r]_\a&\bullet\\
\bullet&&&\bullet\\
\bullet&{}\ar@{}[r]^{\O20}&&\bullet\\
&\bullet\ar@{-}[r]^\a&\bullet}}
\hskip1pt
\boxed{\xymatrix@1@=12pt{
&\bullet\ar@{-}[r]_\a&\bullet\\
\bullet\ar@{-}[d]^\a&&&\bullet\ar@{-}[d]_\a\\
\bullet&{}\ar@{}[r]^{\O40}&&\bullet\\
&\bullet\ar@{-}[r]^\a&\bullet}}$
}
\vtop{\vskip0mm\hbox{\hskip3mm\text{Table 2. Monodromy groups}}
$
\hbox{\boxed{\matrix\\ \phantom{A}\DD_4\phantom{A}\\ \\ \endmatrix}\hskip2pt
\boxed{\matrix\\ \Z_2\oplus\Z_2\\ \\ \endmatrix}\hskip2pt
\boxed{\matrix\\ \phantom{A}\DD_4\phantom{A}\\ \\ \endmatrix} }
$\vskip0mm $
\hbox{\boxed{\matrix\\ \phantom{A}S_3\phantom{A}\\ \\ \endmatrix}\hskip2pt
\boxed{\matrix\\ \phantom{A}\ \Z_2\phantom{A}\ \,\\ \\ \endmatrix}}
$\vskip0mm $
\hbox{\boxed{\matrix\\ \phantom{A}S_4\hskip1pt\phantom{A}\\ \\ \endmatrix}\hskip2pt
\boxed{\matrix\\ \Z_2\oplus\Z_2\\ \\ \endmatrix}\hskip2pt
 \boxed{\matrix\\ \phantom{A}S_4\phantom{A}\\ \\ \endmatrix} }
$
\vskip5mm
\hbox{\hskip19mm
\vtop{\hskip-15mm\hbox{Table 3. Number of orbits}\vskip1mm
\hbox{\boxed{1}\hskip2pt \boxed{2}\hskip2pt\boxed{1} }
\vskip0mm
\hbox{\boxed{2}\hskip2pt\boxed{4}}
\vskip0mm
\hbox{\boxed{1}\hskip2pt
\boxed{2}\hskip2pt \boxed{1} }
}}}}
\endinsert

The edges of graph $\G_\OO$ split into two types:
oval-type labeled by ``$\a$'' on Table 1, and bridge-type labeled by ``$\b$''.
These types indicate
the degeneration of the Hessian curve (collapsing of an oval or a bridge) as the endpoints of an edge are merging.
The numbers of edges of the corresponding types are denoted by $\a(\OO)$ and $\b(\OO)$ and called respectively
the {\it oval-index}
and the {\it bridge-index} of $\OO$.

Our second principal result
(see Theorems \ref{main-theorem} and \ref{coarse-classification})
shows that the pair $\a(\OO),\b(\OO)$ is a complete invariant of coarse deformation equivalence, and
any combination $(\a,\b)\in\{0,2,4\}\times\{0,3,4\}$ except $(4,3)$ is realizable by some Cayley M-octad $\OO$.
%
Our next result is description of the mutual position (adjacency) of the corresponding coarse deformation components $\O{\a}{\b}$:
in Fig.\,1, Table 1 adjacent components stand next to each other (in the same row or column).
After proving chirality of all Cayley octads (see Theorem  \ref{chirality-of-octads}), we deduce adjacency of the pure deformation
components  $\O{\a}{\b}^\pm$.

Finally, in Theorem \ref{real-monodromy} we described the
{\it real monodromy groups} $\Aut(\OO)\subset S_8$ of M-octads $\OO$: they are indicated in Table 2 (the 8 cells
of this table correspond to the cells of Table 1). Here group $\DD_4$ is dihedral of order $8$.
These groups act on the graphs $\Gamma_\OO$ preserving the decoration of edges.
Table 3 gives the number of orbits of $\Aut(\OO)$ acting on $\OO$, whose meaning is clarified in Sect. 6.1.


\section{Preliminaries}

In this Section, we outline some essentials on the Cayley octads and their Hessian curves with spectral theta-charactristics
(for more details see \cite{Dolgachev} and \cite{GH})
and recall a few well-known facts on real plane quartics.

\subsection{Degeneration of Cayley octads and nets of quadrics}
The variety of Cayley octads form a Zariski-open
subset of the Grassmannian of 2-planes in the projective space of quadrics, $\N\subset P(\Sym^2(\C^4))$, namely,
the nets of quadrics $\N$ that have purely zero-dimensional basepoint locus, $\OO=\OO(\N)$.
%
The Hessian curve $\H$ of a Cayley octad $\OO$ is interpreted as the intersection $\H=\N\cap\D$, where $\N$ is the associated
net of quadrics and $\D\subset P(\Sym^2(\C^4))$ the discriminant hypersurface.

In the variety of Cayley octads regular ones form a Zariski-open subset
and singular ones form a certain discriminant locus, whose principal stratum
(of codimension 1) is characterized
in terms of nets $\N$ and their Hessian curves $\H=\N\cap\D$
as follows (cf. \cite{GH}, Lemma 6.4).

\lemma\label{sing-octads}
A singular Cayley octad $\OO$ represents
a non-singular points of the discriminant locus of the variety of octads
if and only if the $\H$ has one node and no other singular points.
Such a node may appear in two ways.
\roster\item
At the point of (simple) tangency of $\N$ with $\D$.
Then the net $\N$ contains one and only one cone with the vertex at some point of $\OO$,
and the Cayley octad $\OO$ has one double point, whereas the other 6 points are ordinary and no 4 of them are coplanar.
\item
At the point of (generic) intersection of $\N$ with the singular locus $\Sing(\D)$ of $\D$.
Then $\N$ contains a reducible quadric (a pair of distinct planes), and
the 8 points of $\OO$ are distinct and admit one and only one splitting into two coplanar quadruples.
\endroster
\endlemma

\proof
The claim on the non-singular points is straightforward, and the description how a node on $\H=\N\cap\D$ may appear is also trivial.
The locus $\Sing(\D)$ parameterizes reducible quadrics, which proves the interpretation formulated in the case (2).
In the remaining case then the description (1) follows.
\endproof

In the case (1), we say that Cayley octad $\OO$ experiences a {\it 2-collision}, and in the case (2)
(two coplanar quadruple of points), a {\it 4-collision}.

The space of real octads is the real locus of the space of complex octads, that is
the corresponding Zariski-open subset in the real Grassmannian of 2-planes in  $P(\Sym^2(\R^4))$.
We will be interested in the part of this space represented by M-octads. Namely, we denote by $\CCC^*$, $\DO2$ and $\DO4$
the strata formed by regular M-octads, by singular ones experiencing 2-collision and 4-collision respectively.
Together they give space $\CCC=\CCC^*\cup\DO2\cup\DO4$.

\lemma\label{interior-walls}
Space $\CCC$ is a manifold with boundary $\DO2$.
The stratum $\DO4$ lies in the interior of $\CCC$ and separates
the connected components of $\CCC^*$.
\endlemma

\proof
As it follows from Lemma \ref{sing-octads},
$\DO2$ as well as $\DO4$ are strata of codimension 1 in the variety of real Cayley octads, so,
it is enough to check what lies from the two sides of them.
The double point of an M-octad $\OO\in\DO2$ can be perturbed into a pair of real points, or a pair of conjugate imaginary ones
by a small real variation of $\OO$, which implies that $\DO2$ is bounded by $\CCC$ from one side.
On the contrary, the number 8 of real points of $\OO\in\DO4$ is preserved after any real variation, so $\OO$ is adjacent to $\DO4$ from both sides.
\endproof

\subsection{The spectral correspondence}
For a non-singular curve $C$, let $\Theta(C)$ be the set of theta-characteristics on $C$.
Consider function
$$h\: \Theta(C)\to\Z/2,\ \ \theta\mapsto \dim H^0(C;\theta)\ \roman{mod}\ 2,\quad \text{ and let}\quad
\Theta_i(C)=h^{-1}(i), i=0,1,$$
which gives a partition $\Theta(C)=\Theta_0(C)\cup\Theta_1(C)$ of theta-characteristics into
{\it even} and {\it odd},
respectively.
Given a regular Cayley octad $\OO$ with the Hessian curve $\H=\N\cap\D$, there is a natural line bundle
$\Cal L_\H\to \H$, whose fiber over $Q_t\in\H\subset\N$ is the null-spaces
of the degenerate quadratic form representing $Q_t$.
This bundle represents linear system $|K_\H+\theta_\OO|$,
where $K_\H$ is the canonical class of $\H$ and $\theta_\OO$ is an even theta-characteristic
called {\it spectral}. This linear system embeds $\H$ to $P^3$ as a sextic, and the image,
$\til\H$, is called the {\it Steinerian curve}.
More precisely, the map $\H\to\til\H$
associates to a singular quadric $Q_t$ (a cone) from the net $\N$ the vertex of $Q_t$,
and so, $\til\H$ lies in the original space $P^3$ containing $\OO$ (see \cite{Dolgachev} for details).
This correspondence $\OO\mapsto(\H,\theta_\OO)$
is bijective up to projective equivalence.

\theorem(Hesse-Dixon)\label{Dixon-th}
For any non-singular plane quartic $C$ and even theta-characteristic $\theta$ on $C$
there exists a unique up to projective transformation of $P^3$ regular Cayley octad $\OO$, such that $C$ is
projectively equivalent to its Hessian curve $C$
and $\theta$ is identified by this equivalence with the spectral theta-characteristic of $\OO$.
\qed\endtheorem

We outline the proof in Sect.\,2.3: it goes back to O.Hesse,
who explained how $\OO$ is reconstructed from $(C,\theta)$,
then A.Dixon \cite{Dixon} elaborated and generalized it by proving a similar correspondence for nets of quadrics in any dimension.
These arguments are applicable also in the real setting
(cf. \cite{DIK}, Th. 3.3.4).
Recall that for a real curve $C$ reality of $\theta\in\Theta(C)$ means that it is preserved invariant under the
involution on $\Theta(C)$ induced by the complex conjugation.

\theorem\label{real-Dixon}
For any real quartic $C$ and an even real theta-characteristic $\theta$,
there exists a unique, up to real projective transformation of $P^3_\R$, real net of quadrics in $P^3_\R$ such that
$C$ is its Hessian curve and $\theta$ its spectral theta-characteristic.
\qed\endtheorem

Pairs $(C,\theta)$, with a non-singular quartic $C\subset P^2$ and an even theta-characteristic $\theta$
will be called {\it even spin quartics}.
The spectral correspondence is not only a bijective map between the projective classes of regular Cayley octads and projective classes of non-singular even spin
quartics, but is an isomorphism of the corresponding moduli spaces, see \cite{Dolgachev}, Remark 6.3.1 and \cite{GH} Sect.\,6.
Restricting it to the real loci we can conclude in particular that their connected components are in one-to-one correspondence.

\corollary\label{coarse-correspondence}
The spectral correspondence induces a bijective correspondence between the coarse deformation classes
of real regular Cayley octads and the deformation classes of real even spin quartics.
\qed\endcorollary

\subsection{Points of a Cayley octad as Aronhold sets for its Hessian quartics
}\label{Aronhold-section}
The line $x_ix_j$ connecting a pair of points of a regular Cayley octad  $\OO=\{x_0,\dots,x_7\}$
determines a line $B_{ij}$ in the plane $\N$ that is a pencil of quadrics whose base-locus
is the union of line $x_ix_j$ with a twisted cubic passing through the other six points of $\OO\sm\{x_i,x_j\}$.
Line $B_{ij}$ is bitangent to the Hessian quartic $\H$ and the divisor formed by the two tangency points
defines an odd theta-characteristic denoted $\theta_{ij}$ (see \cite{Dolgachev}, Theorem 6.3.2).
This gives one-to-one correspondence between the set
$\Theta_1(\H)$ of 28 odd theta-characteristics on $\H$ and the 28 pairs $\{i,j\}\subset\{0,\dots,7\}$.

Any quadruple $\{i,j,k,l\}\subset\{0,\dots,7\}$ defines an even theta-characteristic
$$\theta_{ijkl}=\theta_{ij}+\theta_{ik}+\theta_{il}-K_\H.$$
The choice of a distinguished index $i$ here is not essential, and moreover, the complementary quadruple
$\{0,\dots,7\}\sm\{i,j,k,l\}$ gives the same theta-characteristic as $\theta_{ijkl}$.
This describes $35=\frac12\binom84$ even theta-characteristics of the set $\Theta_1(\H)\sm\{\theta_\OO\}$.
The remaining even theta-characteristic is
$$\theta_\OO=\theta_{01}+\theta_{02}+\dots+\theta_{07}-3K_\H.$$

The set of 7 bitangents $B_{0i}$ corresponding to the
7 summands $\theta_{0i}$ in $\theta_\OO$ form an {\it Aronhold set} (see \cite{Dolgachev}, Sect.6.3.3 for details).
There exist precisely $288=8\times 36$ Aronhold sets representing 36 even theta-characteristics, so that
every chosen $\theta\in\Theta_0(\H)$ is represented precisely by 8 Aronhold sets, and
each Aronhold set match precisely to one vertex of the Cayley octad $\OO$ defined by $(\H,\theta)$, $\theta\in\Theta_0(\H)$.

Namely, given such $(\H,\theta)$, we embed $\H$ to $P^3$ by the
linear system $K_\H+\theta_0$.

Each binangent $B_i$ of a fixed Aronhold set $B_1,\dots, B_7$ representing $\theta$, gives a pair of points on $\til\H$.
The seven lines passing through such pairs of points must be concurrent and intersect at the point of $\OO$ corresponding
to the fixed Aronhold set. This is how $\OO$ is reconstructed from a given even spin quartic.

\subsection{Theta-characteristics as quadratic functions}\label{theta-quadratic}
With $\theta$ we associate
$$q_\theta\:H_1(C;\Z/2)\to\Z/2, \quad x\mapsto h(\theta)+h(\theta+x^*)\mod2,$$
where $x^*\in H^1(C;\Z/2)$ is Poincare dual to $x\in H_1(C;\Z/2)$ (see \cite{Atiyah} and \cite{Mumford}).
It is a quadratic function in the sense that
$$ q_\theta(x+y)=q_\theta(x)+q_\theta(y)+x\cdot y\mod2,
\tag{1}
$$

This gives the well-known identification of the set $\Theta(C)$ with the set of quadratic functions on $H_1(C;\Z/2)$,
compatible with the action of $H^1(C;\Z/2)$ on the both sets.
The function  $h$
is identified then with the
Arf invariant of quadratic functions, $\Arf(q_\theta)\in\Z/2$.
Note that by definition of $q_\theta$
$$
 \Arf(q_{\theta+x^*})=h(\theta+x^*)=h(\theta)+q_\theta(x)=\Arf(q_\theta)+q_\theta(x)\mod2\tag{2}
$$

In the case of Cayley octad $\OO$ and its Hessian quartic $\H$ and $x\in H_1(\H;\Z/2)$ we obtain
$$
q_{\theta_\OO}(x)=\Arf(q_{\theta_\OO+x^*})=h(\theta_\OO+x^*)\mod2\tag{3}
$$

\subsection{Bipartitions of Cayley octads}
For a a regular Cayley octad $\OO=\{x_0,\dots,x_7\}$ we consider its {\it even bipartitions} $\OO=A\cup B$, $A\cap B=\oo$,
into subsets of even cardinalities $|A|$ and $|B|$, and denote by $\LLL(\OO)$
the set formed by 64 such unordered pairs $\{A,B\}$.
 Then the assignments
 $$\{x_i,x_j\}\mapsto \theta_{ij},\quad \{x_i,x_j,x_k,x_l\}\mapsto \theta_{ijkl},\quad \OO\mapsto\theta_\OO$$
 induces a natural one-to-one correspondence $\Phi^{\Theta}\:\LLL(\OO)\to\Theta(\H)$.

Note that $\LLL(\OO)$ is a $\Z/2$-vector space as a subquotient of the power set $\Cal P(\OO)$, with the sum operation
induced by the symmetric difference $\vartriangle$:
$$\{A,B\}+\{C,D\}=\{A\!\vartriangle\! C,A\!\vartriangle\! D\}=\{B\!\vartriangle\! D,
B\!\vartriangle\! C\}$$
and is endowed with the non-degenerate $\Z/2$-valued inner product
 $$\{A,B\}\cdot\{C,D\}=|A\cap B|\mod2.$$

The difference $\theta-\theta'$ for $\theta,\theta'\in\Theta(\H)$ is a 2-torsion element of $\operatorname{Pic}^0(\H)$ and such elements are identified
with the classes in $H^1(\H;\Z/2)$, so, we obtain the associated one-to-one correspondence
$\Phi^H\:\LLL(\OO)\to H^1(\H;\Z/2)$ which is the composition of $\Phi^\Theta$ and the map $\theta\mapsto \theta-\theta_\OO$.
Composing $\Phi^H$ with the Poincare duality we obtain also a map $\Phi_H\:\LLL(\OO)\to H_1(\H;\Z/2)$. 

\proposition\label{identification}
\roster\item
The maps $\Phi^H$ and $\Phi_H$ are linear isomorphisms of $\Z/2$-vector spaces.
\item
Isomorphism $\Phi^\Theta$ identifies $h\: \Theta(\H)\to\Z/2$,
with the map $\LLL(\OO)\to\Z/2$, $\{A,B\}\mapsto \frac12{|A|}\mod2$.
\item
If $x=\Phi_H(\{A,B\})\in H_1(\H;\Z/2)$, then
$q_{\theta\OO}(x)=h(x)=\frac12{|A|}\mod2$.
\item
$\Phi_H$  preserves the inner product, namely, for all $x,y\in H_1(\H;\Z/2)$
$$\Phi_H^{-1}(x)\cdot \Phi_H^{-1}(y)=x\cdot y,\quad \text{ where } x\cdot y \text{ is the homology intersection index}.
$$
\endroster
\endproposition

\proof
(1) Linearity of  $\Phi^H$ (and thus, $\Phi_H$) follows from the relation $\theta_{ij}+\theta_{jk}+\theta_{ik}=\theta_\OO+K_\H$,
or equivalently, $\theta_{ik}=\theta_{ij}+\theta_{jk}-\theta_\OO$
(cf. \cite{GH, Sect.\,6} or \cite{Dolgachev, Theorem 6.3.3}).
For instance, we can deduce from it
$$
\aligned
&\Phi^H(\{x_i,x_j\},\OO\sm\{x_i,x_j\})+\Phi^H(\{x_j,x_k\},\OO\sm\{x_j,x_k\})=\\
&(\theta_{ij}-\theta_\OO)+(\theta_{jk}-\theta_\OO)=
\theta_{ik}-\theta_\OO= \Phi^H(\{x_i,x_k\},\OO\sm\{x_i,x_k\}),
\endaligned
$$
and similarly, linearity in the other cases.
Bijectivity of $\Phi^\Theta$ implies that the linear maps $\Phi^H$ and $\Phi_H$ are isomorphisms.

Part (2) holds by definition of $\Phi^\Theta$, since $\theta_{ij}$ are odd and $\theta_{ijkl}$ are even.

Part (3) follows from relation (3) of Sect. \ref{theta-quadratic}, which taking into account also relation (1) 
applied to $q_{\theta_\OO}$ implies part (4).
\endproof

%

\subsection{Picard-Lefschetz monodromy transformation in $\Theta(\H)$}

\lemma\label{theta-monodromy}
The Picard-Lefschetz monodromy transformation associated with the vanishing class $v\in H_1(\H;\Z/2)$ induces
map $T_v^\Theta\:\Theta(\H)\to\Theta(\H) $
$$\theta\mapsto\theta+(q_\theta(v)+1)v^*.\tag{3}$$
In particular, $\theta$ is preserve invariant if and only if $q_\theta(v)=1$.
\endlemma

\proof
As is well-known, the monodromy action in $H_1(\H;\Z/2)$
associated with $v$ is the Picard-Lefschetz
transformation
$x\mapsto x+(x\cdot v)v$,
so, the corresponding monodromy transforms in $\Theta(\H)$,
$q_\theta\mapsto q'_\theta$ should satisfy modulo 2 relation
\newline
$q'_\theta(x)=q_\theta(x+(x\cdot v)v)=q_\theta(x)+(x\cdot v)q_\theta(v)+(x\cdot v)^2=
q_\theta(x)+(q_\theta(v)+1)v^*(x).$
\endproof

\proposition\label{octadic-monodromy}
The Picard-Lefschetz transformation $T^\Theta_v$ preserves the pairity of theta-characteristics $\theta$ (Arf-invariant of $q_\theta$).
The action induced by $T^\Theta_v$ in $\LLL(X)$ via $\Phi^\Theta$ depends on the bipartition
$\Phi^{-1}_H(v)=\{A_v,B_v\}\in\LLL(X)$ as follows.
\roster\item
If one of the sets $A_v$ or $B_v$ is a 2-elements set
$\{x_i,x_j\}\subset\OO$ then the action of $T^\Theta_v$ in $\LLL(X)$ is induced by the transposition of $x_i$ and $x_j$.
\item
If $A_v$ and $B_v$ are 4-element sets, then $T^\Theta_v$ replaces 2-element subsets $\{x_i,x_j\}\subset A_v$ 
$($or $\{x_i,x_j\}\subset B_v$\,$)$
by their complements $A_x\sm\{x_i,x_j\}$ $($respectively $B_v\sm\{x_i,x_j\}$\,$)$ and keep $\{x_i,x_j\}$ unchanged if it has one common point with $A_v$ and $B_v$.
\endroster
\endproposition

\proof
Applying formula (3) and
using relation (1), (2) we obtain
$$\Arf(q'_\theta)=\Arf(q_\theta)+q_\theta((q_\theta(v)+1)v)=\Arf(q_\theta)\!\mod2,$$
since $q_\theta((q_\theta(v)+1)v)=(q_\theta(v)+1)q_\theta(v)$,
which means preserving of pairity of $\theta$.

The action on $\{A,B\}\in\LLL(\OO)$ induced by $T^\Theta_v$ due to Lemma \ref{theta-monodromy},
in terms of $\{A_v,B_v\}$ is as follows:
$$
\{A,B\}\mapsto \{A,B\}+\e\{A_v,B_v\},\quad \text{where } \e=\frac{|A_v|}2+|A\cap A_v|+1\mod2,
$$
which gives $\e=|A\cap A_v|\mod2$ if $A_v=\{x_i,x_j\}$, and so, in the case (1) $T^\Theta_v$ is induced by transposition of $x_i$ and $x_j$.
In the case (2) $\e=|A\cap A_v|+1\mod2$ and thus, $T^\Theta_v$ acts also as is described.
\endproof

Transformations part (2) of Proposition \ref{octadic-monodromy} were named {\it bifid substitutions} by Cayley, who found them as
the monodromy of 4-collisions in the corollary below.

\corollary\label{collision-correspondence}
Let $\H$ be the Hessian quartic of a regular Cayley octad $\OO$ with the spectral theta-characteristic $\theta\in\Theta(\H)$ and
associated quadratic function $q_\theta$.
Assume that $v\in H_1(\H;\Z/2)$ is a vanishing class associated with some nodal degeneration of $\H$, $v=\Phi_H(\{A_v,B_v\})$. Then
$q_\theta(v)=\frac{|A|}2\mod2$ and
the degeneration of $\OO$ corresponding to the nodal degeneration of $\H$ is
\roster\item
a 2-collision of $x_i,x_j\in\OO$ that form 2-element set, $A_v$ or $B_v$, in the case of $q_\theta(v)=1$;
\item
a 4-collision involving $A_v$ and $B_v$, in the case of $q_\theta(v)=0$.
\endroster
\endcorollary

\proof
It follows from the description of the Picard-Lefschetz transformation corresponding
to the vanishing class $v$ in terms of $\OO$ given in Proposition \ref{octadic-monodromy}.
\endproof

\subsection{Real non-singular quartics}
The real deformation classification of non-singular real plane quartics $C$ goes back to F.Klein:
there exist 6 real deformation classes characterized by the number $0\le k\le 4$ of
components of $C_\R$ called {\it ovals}
and in the case $k=2$ by their mutual position:
in one class
two ovals bound disjoint pair of discs in $P^2_\R$, and in the other two ovals are {\it nested}:
one oval bounds a disc containing another oval.
By {\it M-quartics} we mean real non-singular quartics with the maximal number 4 of ovals.
The ovals of such quartics can be arbitrarily permuted by the monodromy.

\proposition{\bf(Klein)}\label{Klein}
M-quartics form one real deformation class. The ovals of an M-quartic can be arbitrarily permuted by the
real deformation monodromy: the corresponding monodromy group is the symmetric group $S_4$.
\qed\endproposition

\proposition\label{discrepancy-correspondence}
A real regular Cayley octad $\OO$ has $8-2d$ real points, $0\le d\le 3$,
if and only if its Hessian quartic $\H$ has $4-d$ real ovals, which (for $d=2$) are not nested.
\endproposition

\proof
The real structure from $\H$ can be lifted to
the del Pezzo surface $Z$ obtained by double covering of $P^2$ branched along $\H$, so that
the real locus $Z_\R$ is projected to the non-orientable region of $P^2_\R$ bounded by $\H_\R$.
Then $\chi(Z_\H)=2(\chi(P^2_\R)-k)=2-2k$, where $k$ is the number of components of $\H_\R$.
On the other hand, if $\OO$ has $8-2d>0$ real points, then in the blowup model of $Z$, its real  locus
$Z_\R$ is obtained by blowing up $P^2_\R$ at $8-2d-1$ points, thus, $\chi(Z_\R)=1-(8-2d-1)=2-2(4-d)$.
Thus, $k=4-d$.
\endproof

\section{Real Hessian curve with real $\theta$-characteristics}

\subsection{The ovals and bridges of M-quartics}
Given an M-quartic $C\subset P^2$ we consider {\it real vanishing cycles},
which are $\conj$-invariant simple closed curves on $C$ that can be contracted to the nodal point by some real
degeneration of $C$.
Such cycles include
the the four ovals, $\T_0,\dots, \T_3\subset C_\R$:
contraction of an oval yields a node of {\it solitary} type.

Another kind of vanishing cycles are {\it bridges} $b_{ij}$, $0\le i<j\le3$, connecting ovals $\T_i$ and $\T_j$.
Namely, the quartic splitting into four real lines in general position in $P^2$ can be
perturbed into M-quartic, and
the six intersection points of the lines yield the six bridges $b_{ij}$,
which are real vanishing cycles having precisely two real points at the intersection with ovals  $\T_i$ and $\T_j$.
Connectedness of the space of M-quartics yields such bridges for any given M-quartic $C$ and Proposition \ref{oval-bridge-bases}
below shows that the classes $[b_{ij}]\in H_1(C;\Z/2)$ are well-defined (unique).


Consider $\pm1$-eigengroups $H_1^\pm(C)=\{x\in H_1(C)\,|\,c(x)=\pm x\}$
along with their images  $H_1^\pm(C;\Z/2)\subset H_1(C;\Z/2)$ under the modulo 2 reduction homomorphism.
It is a well-known fact that the action of $c$ in $H_1(C;\Z/2)$ is trivial for M-curves, which leads to the
direct sum decompositions
$$H_1(C)=H_1^+(C)\oplus H_1^-(C),\quad
 H_1(C;\Z/2)=H_1^+(C;\Z/2)\oplus H_1^-(C;\Z/2).$$

\proposition\label{oval-bridge-bases}
(1) Classes $[\T_i]$, $0\le i\le3$ span $\Z/2$ vector space
$H_1^+(C;\Z/2)$ with the only relation $[\T_0]+\dots+[\T_3]=0$.
In particular, any three of these classes form a basis.

(2)
The homology classes $[b_{ij}]\in H_1(C;\Z/2)$, $0\le i<j\le3$ span  $H_1^-(C;\Z/2)$ and satisfy relations
 $[b_{ij}]+[b_{jk}]+[b_{ik}]=0$
 for any triple of indices $0\le i<j<k\le3$
 In particular,
classes $[b_{0i}]$, $1\le i\le3$ form a basis of $H_1^-(C;\Z/2)$.

(3) Any real nodal degeneration which leads to merging of ovals $\T_i$ and $\T_j$ at a cross-like node
is described by the vanishing class $[b_{ij}]\in H_1(C;\Z/2)$.
 In particular, bridge classes
are independent of the choice of a degeneration of an M-quartic $C$ into four real lines.
\endproposition

\proof
If $C$ is an M-quartic, then the complement $C\sm C_\R$ splits into
two connected components, $C_1,C_2$, topological spheres with 4 holes permuted by the complex conjugation,
which yield the homology relation in (1) and trivially implies that classes $[\T_i]$ span $H_1^+(C;\Z/2)$.

Since bridge $b_{ij}$ is the union of a path connecting $\T_i$ with $\T_j$ in $C_1$ with the $\conj$-symmetric path in $C_2$,
we can also trivially deduce (2).

Finally, $[b_{ij}]$ is determined by the modulo $2$ intersection indices $b_{ij}\cdot\T_i=b_{ij}\cdot\T_j=1$, $b_{ij}\cdot\T_k=0$, for $k\ne i,j$,
which gives (3).
\endproof

\subsection{The real monodromy action on theta-characteristics}
Consider the variety $\Cal C$ of spin M-quartics, $(C,\theta)$, $\theta\in\Theta(C)$. It splits
into two components, $\Cal C=\Cal C^{ev}\cup\Cal C^{odd}$, where for $C^{ev}$ and $\Cal C^{odd}$
$\theta$ is respectively even and odd.
Each of these two components is furthermore partitioned as
$$\Cal C^{ev}=\bigcup_{\a,\b}\Cal C^{ev}_{\a,\b},\quad \Cal C^{odd}=\bigcup_{\a,\b}\Cal C^{odd}_{\a,\b},$$
where subscript $\a$ in $\Cal C^{ev}_{\a,\b}$ and $\Cal C^{odd}_{\a,\b}$
stands for the numbers of ovals $\T_i$ for which $q_\theta([\T_i])=1$, while $\b$
counts the number of bridge classes $[b_{ij}]$ with $q_\theta([b_ij])=1$.

In this Subsection we will prove the following result.

\theorem\label{main-theorem}
Space $\Cal C$ has precisely 11 connected components: 8 components in $\Cal C^{ev}$ and 3 components in $\Cal C^{odd}$.
The components of $\Cal C^{ev}$ are $\Cal C^{ev}_{\a,\b}$ where $a\in\{0,2,4\}$, $b\in\{0,3,4\}$, with exception of the pair $(\a,\b)=(2,4)$,
and the components of $\Cal C^{odd}$ are $\Cal C^{odd}_{\a,\b}$ where $(\a,\b)$ is $(2,4)$, $(2,3)$, or $(4,3)$.
\endtheorem

We start proving it with two lemmas.

%
%

%

\lemma\label{64theta}
The 64 theta-characterisitics $\theta\in\Theta(C)$ on an M-quartic $C$ are classified by the three $\Z/2$-values $q_\theta([\T_i])$
and three $\Z/2$-values $q_\theta([b_{0i}])$, $i=1,2,3$.
The Arf-invariant of $q_\theta$ is equal to
$$Arf(q_\theta)=\sum_{1=3}^3q_\theta([\T_i])q_\theta([b_{0i}]).$$
\endlemma

\proof
It follows trivially from Proposition \ref{oval-bridge-bases} and from the definition of the Arf-invariant.
\endproof

So, we can encode the theta-characteristics  $\theta\in\Theta(C)$ by binary $2\times3$-matrices
$\left[\matrix
q_\theta([\T_1])&q_\theta([\T_2])&q_\theta([\T_3])\\
q_\theta([b_{01}])&q_\theta([b_{02}])&q_\theta([b_{03}])
\endmatrix\right]$.
Symmetric group $S_4$ permuting the ovals (or in the other words, changing ordering $\T_0,\dots,\T_3$)
acts on the set $\Theta(C)$. Namely, $\sigma\in S_4$
induces an action on the bridges classes sending $[b_{ij}]$ to  $[b_{\sigma(i)\sigma(j)}]$,
which determines the induced action on the set of 64 binary $2\times3$-matrices.
We will list the orbits of this actions.

\lemma\label{codes-of-theta-orbits}
There exist precisely 11 orbits of the $S_4$-action
on $\Theta(C)$: eight in $\Theta_0(C)$ and three in $\Theta_1(C)$.
In terms of binary matrices they are represented by
\roster\item
eight matrices for $\theta\in\Theta_0(C)$

\vtop{\hbox{\hskip-1pt
$\left[\matrix
0&0&0\\
1&0&1
\endmatrix\right]$,
$\left[\matrix
0&1&0\\
1&0&1
\endmatrix\right]$,
$\left[\matrix
1&1&1\\
1&0&1
\endmatrix\right]$,}\vskip1mm
\hbox{$\left[\matrix
0&0&0\\
0&0&1
\endmatrix\right]$,
$\left[\matrix
0&1&0\\
0&0&1
\endmatrix\right]$,}\vskip1mm
\hbox{$\left[\matrix
0&0&0\\
0&0&0
\endmatrix\right]$,
$\left[\matrix
0&1&0\\
0&0&0
\endmatrix\right]$,
$\left[\matrix
1&1&1\\
0&0&0
\endmatrix\right]$,}}

\item
three matrices for $\theta\in\Theta_1(\C)$

\hbox{
$\left[\matrix
0&1&1\\
1&0&1
\endmatrix\right]$,
$\left[\matrix
0&1&1\\
0&0&1
\endmatrix\right]$,
$\left[\matrix
1&1&1\\
0&0&1
\endmatrix\right]$.}
\endroster
\endlemma

\proof
Enumeration of orbits is a straightforward exercise.
The pairity of $\theta\in\Theta(C)$ represented by the matrices listed are determined by Lemma \ref{64theta}.
\endproof

\demo{Proof of Theorem \ref{main-theorem}}
Connectedness of the space of M-quartics implies immediately that
the set of connected components of the space $\Cal C$
is in one-to-one correspondence with the set of orbits of the real deformation monodromy action
on the set $\Theta(C)$ for any chosen M-quartic $C$.
Since the real monodromy group acting on the ovals of $C$ is $S_4$ (see \ref{Klein}),
Lemma \ref{codes-of-theta-orbits} implies the first part of the Theorem.
For the second part we need just to determine the values of $\a$ and $\b$
(numbers of non-zero values of $q_\theta$ on the 4 ovals and 6 bridges)
for the 11 matrices in Lemma \ref{codes-of-theta-orbits}.
It is straightforward because
both subgroups $H^\pm(C;\Z/2)$ are isotropic with respect to the intersection form, and thus,
the restriction of $q_\theta$ to each of them is linear.
Namely, in the list (1) values of $\a$
are respectively $0$, $2$ and $4$ in the first, the second and the third columns, while
the values of $\b$ are respectively $4$, $3$ and $0$ in the first, second and third rows.

In part (2) the three matrices represent pairs $(\a,\b)$ which are, respectively, $(2,4)$, $(2,3)$, and $(4,3)$.
\qed\enddemo

\subsection{Theta-diagrams of M-quartics}
We give below a simple graphical interpretation of the $2\times3$-matrices listed in Lemma \ref{64theta}, which
reveals the geometry behind the enumeration of orbits in Lemma \ref{codes-of-theta-orbits} and simplifies it.

\midinsert
\topcaption{Fig.\,2. Ovals and bridges of M-quartics}\endcaption
\vskip-2mm
\hbox{\hskip-5mm\vtop{\hsize85mm\hbox{\hskip40mm Even $\theta$}
\hbox{\vtop{\hsize13mm\vskip6mm $\,\,\beta=4$}
$\vtop{\hsize10mm\boxed{\xymatrix@1@=10pt{
{}&&&{}\\
&\bullet\ar@{-}[ul]\ar@{.}[r]\ar@{.}[d]&\bullet\ar@{.}[d]\ar@{-}[ur]\\
&\bullet\ar@{-}[dl]\ar@{.}[r]&\bullet\ar@{-}[dr]\\
{}&&&{}
}}}\hskip3pt
\vtop{\hsize10mm\boxed{\xymatrix@1@=10pt{
{}&&&{}\\
&\bullet\ar@{-}[ul]\ar@{.}[r]\ar@{.}[d]&\circ\ar@{.}[d]\ar@{-}[ur]\\
&\circ\ar@{-}[dl]\ar@{.}[r]&\bullet\ar@{-}[dr]\\
{}&&&{}
}}}\hskip3pt
\vtop{\hsize10mm\boxed{\xymatrix@1@=10pt{
{}&&&{}\\
&\circ\ar@{-}[ul]\ar@{.}[r]\ar@{.}[d]&\circ\ar@{.}[d]\ar@{-}[ur]\\
&\circ\ar@{-}[dl]\ar@{.}[r]&\circ\ar@{-}[dr]\\
{}&&&{}
}}
}$}
\vskip1pt
\hbox{
\vtop{\hsize12mm \vskip6mm$\b=3$}
$\vtop{\hsize10mm\boxed{\xymatrix@1@=10pt{
{}&&&{}\\
&\bullet\ar@{.}[ul]\ar@{-}[r]\ar@{-}[d]&\bullet\ar@{.}[d]\ar@{-}[ur]\\
&\bullet\ar@{-}[dl]\ar@{.}[r]&\bullet\ar@{.}[dr]\\
{}&&&{}
}}}\hskip3pt
\vtop{\hsize10mm\boxed{\xymatrix@1@=10pt{
{}&&&{}\\
&\bullet\ar@{.}[ul]\ar@{-}[r]\ar@{-}[d]&\circ\ar@{.}[d]\ar@{-}[ur]\\
&\circ\ar@{-}[dl]\ar@{.}[r]&\bullet\ar@{.}[dr]\\
{}&&&{}
}}}
$}
\vskip1pt
\hbox{
\vtop{\hsize12mm \vskip6mm$\b=0$}
$\vtop{\hsize10mm\boxed{\xymatrix@1@=10pt{
{}&&&{}\\
&\bullet\ar@{-}[ul]\ar@{-}[r]\ar@{-}[d]&\bullet\ar@{-}[d]\ar@{-}[ur]\\
&\bullet\ar@{-}[dl]\ar@{-}[r]&\bullet\ar@{-}[dr]\\
{}&&&{}
}}}\hskip3pt
\vtop{\hsize10mm\boxed{\xymatrix@1@=10pt{
{}&&&{}\\
&\bullet\ar@{-}[ul]\ar@{-}[r]\ar@{-}[d]&\circ\ar@{-}[d]\ar@{-}[ur]\\
&\circ\ar@{-}[dl]\ar@{-}[r]&\bullet\ar@{-}[dr]\\
{}&&&{}
}}}\hskip3pt
\boxed{\xymatrix@1@=10pt{
{}&&&{}\\
&\circ\ar@{-}[ul]\ar@{-}[r]\ar@{-}[d]&\circ\ar@{-}[d]\ar@{-}[ur]\\
&\circ\ar@{-}[dl]\ar@{-}[r]&\circ\ar@{-}[dr]\\
{}&&&{}
}}$}
\vskip1mm
\hbox{\hskip20mm$\a=0$\hskip12mm $\a=2$\hskip12mm$\a=4$}
}
\hskip5mm
\vtop{\hbox{\hskip5mm {Odd }$\theta$}
\hbox{$\boxed{\xymatrix@1@=10pt{
{}&&&{}\\
&\bullet\ar@{-}[ul]\ar@{.}[r]\ar@{.}[d]&\bullet\ar@{.}[d]\ar@{-}[ur]\\
&\circ\ar@{-}[dl]\ar@{.}[r]&\circ\ar@{-}[dr]\\
{}&&&{}
}}$}
\vskip1pt
\hbox{$\boxed{\xymatrix@1@=10pt{
{}&&&{}\\
&\bullet\ar@{.}[ul]\ar@{-}[r]\ar@{-}[d]&\bullet\ar@{.}[d]\ar@{-}[ur]\\
&\circ\ar@{-}[dl]\ar@{.}[r]&\circ\ar@{.}[dr]\\
{}&&&{}
}}\hskip2pt
\boxed{\xymatrix@1@=10pt{
{}&&&{}\\
&\circ\ar@{.}[ul]\ar@{-}[r]\ar@{-}[d]&\circ\ar@{.}[d]\ar@{-}[ur]\\
&\circ\ar@{-}[dl]\ar@{.}[r]&\circ\ar@{.}[dr]\\
{}&&&{}
}}
$}\vskip2mm
\hbox{$\bullet\ \text{ black oval: }\quad\quad q_\theta([\T_i])=0$}
\hbox{$\circ\ \text{ white oval: }\quad\quad q_\theta([\T_i])=1$}
\hbox{\xymatrix@1@=15pt{
{}\ar@{-}[r]&{}\ \text{ black bridge:\ \ } q_\theta([b_{ij}])=0}}
\hbox{\xymatrix@1@=15pt{
{}\ar@{.}[r]&{}\ \text{ white bridge:\ \ }  q_\theta([b_{ij}])=1}}
}
}

\endinsert

The four ovals $\T_i$ of an M-quartic $C$ are represented by
four vertices depicted as small circles on the plane and the six bridges are represented as edges depicted as line segments
connecting these circles pairwise: $b_{ij}$ connects circle $\T_i$ with $\T_j$ in $\Rp2$, so that two of the six bridges pass ``through infinity'',
as it is shown on Fig.\,2.
 One can view this diagram as a complete graph $K_4$ embedded in $\Rp2$.

A theta-characteristic $\theta\in\Theta(C)$ can be described as a
coloring of such diagram: black color for
an oval $\T_i$ or bridge $b_{ij}$ means that
$q_\theta$ takes value $0$ on the corresponding class $[\T_i]$ or $[b_{ij}]$,
and white means values $1$. Black colored ovals are depicted as filled circles,
and black bridges are usual edges, while white ovals and bridges look like empty circles
and dotted edges.
A diagram decorated in this way will be called the {\it theta-diagram} of $\theta\in\Theta(C)$.
Such diagram with even $\theta$ will be denoted $D_{\a\b}$, where $\a$ and $\b$ are respectively the numbers of white vertices (ovals)
and white edges (bridges) in it.


\subsection{Coarse Deformation Components}

\theorem\label{coarse-classification}
There exist eight coarse deformation classes of regular Cayley M-octads
which are in spectral correspondence with the eight
components $\Cal C^{ev}_{\a\b}$ of $\Cal C^{ev}$ listed in Theorem \ref{main-theorem}.
\endtheorem

\proof
Theorem \ref{real-Dixon} with Corollary \ref{coarse-correspondence} and
Proposition \ref{discrepancy-correspondence} (in the case of $d=0$) give a bijective correspondence between
the set of coarse deformation classes
of real regular M-octads and the set of deformation classes of
non-singular real even spin M-quartics. So, it is left to apply Theorem  \ref{main-theorem}
 that enumerates them.
\endproof

Let us denote by $\CCC_{\a\b}$ the coarse deformation class corresponding to $\Cal C^{ev}_{\a\b}$ by Theorem \ref{coarse-classification}.

\subsection{Adjacency of the coarse deformation components}
Two real deformation classes of M-octads (connected components of $\CCC^*$)
lying on the opposite sides from a wall of $\DO4$
(formed by M-octad with a 4-collision) are said to be {\it adjacent}.
A pair of coarse deformation classes $\CCC_{\a\b}$
containing such an adjacent pair of components will be called adjacent too.

Adjacency of classes  $\CCC_{\a\b}$ are analyzed below through adjacency of the corresponding components $\Cal C^{ev}_{\a\b}$.
Namely, let us fix an M-quartic $C$ with the ovals $\T_i$, $i=0,\dots,3$ and bridges
$b_{ij}$, $0\le i<j\le3$.
We say that theta-characteristics $\theta$ and $\theta'$ in $\Theta_0(C)$
are {\it adjacent theta-characteristics} if
the pairs $(C,\theta)$ and $(C,\theta')$
 represent adjacent components   $\Cal C^{ev}_{\a\b}$ and $\Cal C^{ev}_{\a'\b'}$
(corresponding to adjacent classes  $\CCC_{\a\b}$ and  $\CCC_{\a'\b'}$).

\proposition\label{adjacent-theta}
\roster\item
Group $H_1(C;\Z/2)$ contains precisely 10 real vanishing classes: four oval-classes $[\T_i]$
and six bridge-classes $[b_{ij}]$.
\item
A pair of theta-characteristics $\theta,\theta'\in\Theta_0(C)$ are adjacent if and only if
$\theta'=\theta+v^*$, where $v^*$ is Poincare dual to some real vanishing class $v\in H_1(C;\Z/2)$
satisfying the condition $q_{\theta}(v)=0$.
\endroster\endproposition

\proof
Any nodal degeneration of $C$ is either a contraction of some oval  $\T_i$ with the vanishing class $[\T_i]$,
or merging of two ovals, $\T_i$ and $\T_j$ with the vanishing class
 $[b_{ij}]$ (see Proposition \ref{oval-bridge-bases}(3)). This proves part (1).

By Lemma \ref{interior-walls} and Corollary \ref{collision-correspondence},
the condition $q_{\theta}(v)=0$ in part (2) means that the corresponding to this nodal degeneration wall
is internal (deformation components of spin M-quartics lie on the both sides of it).

For proving the wall-crossing formula $\theta'=\theta+v^*$, recall that
the forgetful map
$(C,\theta)\mapsto C$ from the variety of spin quartics to the space of quartics
is unramified near nodal quartic if for the vanishing class $v$ we have  $q_\theta(v)=1$
and has ramification of index 2 along the stratum of nodal quartics with $q_\theta(v)=0$.
It follows for instance from Lemma \ref{theta-monodromy}.
The Picard-Lefschetz transformation in $\Theta(C)$, as it was observed in Lemma \ref{theta-monodromy}, is
$q_\theta\mapsto q_\theta+(q_\theta(v)+1)v^*$, which in the case of $q_\theta(v)=0$ is just just adding $v^*$,
so, it implies $q_{\theta'}(v)=q_{\theta}(v)+v^*(v)=q_{\theta}(v)$.

Finally, it is left to notice that a loop around a wall in the space of spin quartics
after lifting to the covering space with the ramification of index 2
becomes a path into the adjacent deformation component
and the Picard-Lefschetz formula becomes the wall-crossing formula for
the markings that define a given branched covering.
\endproof

In terms of M-quartic theta-diagrams $D_{\a\b}$ and $D_{\a'\b'}$ of adjacent deformation classes
this proposition means that  $D_{\a'\b'}$ is obtained from  $D_{\a\b}$ after one of the following two modifications (see Fig.3):
\midinsert
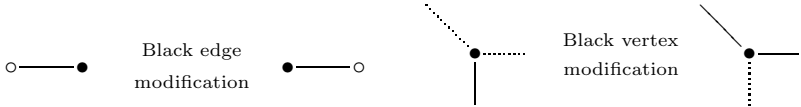
\topcaption{Fig. 3. Black edge and vertex modifications}\endcaption
\vskip-3mm
$\hbox{
\xymatrix@1@=20pt{\\
\circ\ar@{-}[r]&\bullet\ &{}\ar@{}^{\text{Black edge}}_{\text{modification}}[r]{}\ &&{}\bullet\ar@{-}[r]&\circ}}\quad\quad
\hbox{\xymatrix@1@=15pt{
{}\\
&\bullet\ar@{.}[r]\ar@{-}[d]\ar@{.}[ul]&{}\quad
&{}\ar@{}^{\text{Black vertex}}_{\text{modification}}[r]&{}\\
&{} }
\ \xymatrix@1@=15pt{
{}\\
&\bullet\ar@{-}[r]
\ar@{.}[d]\ar@{-}[ul]&{}\\
&{} }}
$
\endinsert
\roster\item
{\it Black edge modification:}
each endpoint, $\T_i$ and $\T_j$, of some black edge $\b_{ij}$ changes its color, while the other 2 vertices and all 6 edges preserve.
\item{\it Black vertex modification:}
each edge $b_{ij}$ incident to some black vertex $\T_i$ changes its color, while the other 3 edges and all 4 vertices preserve.
\endroster

The following theorem shows that adjacency of coarse deformation classes  $\O{\a}{\b}$
corresponds to vertical and horizontal adjacency of theta-diagrams $D_{\a\b}$ on Fig.\,2.

\theorem\label{adjacenct-theta-diagrams}
Two different coarse deformation classes of M-octads $\O{\a}{\b}$ and $\O{\a'}{\b'}$
are adjacent if and only if one of the following holds
\roster
\item $\a=\a'$ and numbers in the pair $\{\b,\b'\}\subset\{0,3,4\}$ have odd difference,
\item
 $\b=\b'$ and $|\a-\a'|=2$.
\endroster

Besides, a coarse deformation component  $\O{\a}{\b}$ is adjacent to itself if and only if $(\a,\b)$ is either $(2,3)$ or $(2,0)$.
\endtheorem

\proof
If the vanishing cycle $v$ representing the wall between adjacent deformation classes is a bridge-class $b_{ij}$, then
according to Lemma \ref{adjacent-theta}
the value of $q_\theta$ changes on the two oval-classes $\T_i$ and $\T_j$ and is preserved on the other oval classes and on all
bridge-classes. Then $\b=\b'$ and $|\a-\a'|$ is either $2$ or $0$. Moreover, $0$ may appear only if some black edge has endpoints of different colour,
which is possible only for $\O{2}{3}$ and $\O{2}{0}$.

If $v$ is an oval-class $\T_i$, then for the same reason $q_\theta$ changes
on the three bridge-classes $[b_{ij}]$, $0\le j\le 3$, $j\ne i$ and so, $\a=\a'$ and $|\b-\b'|$ is either $1$ or $3$.
\endproof

\corollary\label{total-adjacency}
Coarse deformation classes of regular M-octads are related through a finite number of wall-crossings
(passing to an adjacent coarse deformation class of M-octads).
\qed\endcorollary

\subsection{Graphs $\G_\OO$}
Assume that $\OO$ is a regular M-octad and $(\H,\theta)$ is the associated even spin Hessian quartic.
Consider a pair of real vanishing cycles
$v_1,v_2\in H_1(\H;\Z/2)$ representing
2-collisions of $\OO$, and let $e_1$, $e_2$ be the corresponding edges of $\G_\OO$.

\lemma\label{matching} $(1)$
Edges $e_1$ and $e_2$ are different for different $v_1$ and $v_2$.

$(2)$ The intersection index $v_1\centerdot v_2$ is $1$ if and only if $e_1$ and $e_2$ have a common vertex;
in particular, edges $e_1$ and $e_2$ are disjoint if $v_1$, $v_2$ are
both oval-cycles, or both bridge-cycles.
\endlemma


\demo{Proof of Lemma \ref{matching} }
It follows immediately from that $\Phi_H$ preserves $\Z/2$-inner product by
Proposition \ref{identification}.
\qed\enddemo

\proposition\label{Collision_graphs_proposition}
The graphs $\G_\OO$ for $\OO\in\O{\a}{\b}$ have combinatorial types as presented in Table 1 of Fig.\, 1.
\endproposition

\proof
By Proposition \ref{identification}, white ovals and bridges of $\H$ represent edges of graph $\G_\OO$,
so that an incident pair of oval end bridge represents a pair of adjacent edges.
Note that white ovals and edges shown on each of the theta-diagrams of Fig.\,2 split into several chains
formed by consecutively incident ovals and bridges which alternate in the chain.
The corresponding edges in $\G_\OO$ then split into chains of consecutively adjacent edges.
Applying it to each of the eight theta-diagrams, we obtain the graphs $\G_X$ presented in Table 1 of Fig.\,1.

For instance, in the (least trivial) case of $\OO\in\O{4}{4}$, there is one chain forming a cycle,
$\T_0,b_{01},\T_1,b_{12},\dots,\T_7,b_{07}$ (if the ovals are numerated cyclically). This cycle should represent
an cycle on $\G_X$ forming an octagon shown in the top-right cell in Table 1 of Fig.\,1.
Analysis of the other cases is similar.
\endproof

\section{Pure deformation classification}

\subsection{Chirality of configurations in $\Rp3$}
By a {\it simple n-configuration} in $\Rp3$, $n\ge4$, we mean
an $n$-point subset $A\subset\Rp3$
not containing coplanar quadruples of points. A pair of such configurations are said to be {\it deformation equivalent} if they
can be connected by a continuous family of simple $n$-configurations.
A {\it mirror partner} of $A$ is a configuration $\bar A$ obtained from $A$
by an orientation reversing projective transformation (for instance, reflection across a plane).
If a simple $n$-configuration $A$ is deformation equivalent to $\bar A$, it is called
{\it achiral} and otherwise {\it chiral}.
The following observation belongs to O.Viro and V.Kharlamov (see \cite{Viro2}).

\proposition\label{6-7-chirality}
All simple configurations of 6 and 7 points in $\Rp3$ are chiral.
\endproposition

\proof
Following \cite{Viro2}, with a triple of skew lines $\{\ell_1,\ell_2,\ell_3\}$ in $\Rp3$ we associate sign
$$\lk(\ell_1,\ell_2,\ell_3)=\lk(\bar\ell_1,\bar\ell_2)\lk(\bar\ell_1,\bar\ell_3)\lk(\bar\ell_2,\bar\ell_3)\in\{+1,-1\}$$
where $\bar\ell_i$ is the line $\ell_i$ endowed with an arbitrary orientation and
$\lk(\bar\ell_i,\bar\ell_j)=\pm1$ is the normalized linking number in $\Rp3$ (normalization means multiplication by 2, since
the usual linking number is $\pm\frac12$).
 This triple index is clearly independent of the choice of the order of lines $\ell_1$,\dots,$\ell_3$ and of their orientations, but alternates
if we change the orientation of $\Rp3$,
or reflect a triple of lines across a plane.

Given a simple 6-configuration $A$, we consider a triple of skew lines connecting the points of $A$ pairwise.
We can obtain 15 such triples of lines corresponding to 15 splitting of 6 points into 3 pairs, and the product
of the corresponding 15 signs $\lk(\ell_1,\ell_2,\ell_3)$ is a deformation invariant of simple 6-configuration $A$, denoted $\sign_6(A)$.
Then $\sign_6(A')=-\sign_6(A)$ if $A'$ is a mirror of $A$, and so, $A$ is chiral.

For a simple 7-configuration $B$ we obtain seven signs for its 6-subconfigurations, and the product of these signs, $\sign_7(B)\in\{+1,-1\}$,
is a deformation invariant,
such that $\sign_7(B')=-\sign_7(B)$ for a mirror image $B'$ of $B$.
\endproof

Each point $x$ of a simple 8-configuration $\OO$ can be also equipped with a sign, $\sign(X,x)=\sign_7(X\sm\{x\})$.
Then, if $\OO'$ is the mirror partner of $\OO$ and $x'\in\OO'$ the vertex corresponding to $x$, we have $\sign(\OO',x')=-\sign(\OO,x)$.

\theorem\label{chirality-of-octads}
All points of a regular Cayley M-octads are equipped with the same sign.
In particular, any regular Cayley M-octad is chiral. Thus,
there exist precisely 16 pure deformation classes of them.
\endtheorem

We complete its prove towards the end of section 4.

\subsection{Plane 4-moves for 6 and 7-configurations}
We say that two deformation classes of simple $n$-configurations, $n\ge4$,
are {\it adjacent} if they lie from the opposite sides of the codimension 1 stratum formed by
$n$-configurations with a coplanar quadruple of points (a deformation class can be adjacent to itself
if this stratum is one-sided or the same deformation class lie from the both sides).

\lemma\label{4-collision}
Assume that simple $n$-configurations $A_0$ and $A_1$ represent adjacent deformation classes.
Then $\sign_n(A_0)=-\sign_n(A_1)$ if $n=6$ or $n=7$.
In particular, for such $n$ self-adjacent deformation classes do not exist.
\endlemma

\proof
Since $A_0$ and $A_1$ are adjacent, they can be connected
with a path $A_t=\{x_1(t),\dots,x_n(t)\}$, $t\in[0,1]$, in which configurations $A_t$ are simple for all $t$ except one
value $t=t_0$, and for this value
precisely one quadruple of points $x_i(t)$ become coplanar.
 These four points can be connected pairwise by two lines in three ways.
The linking number between the corresponding pairs of lines (with some auxiliary orientations) for $t\ne t_0$
alternates
as we cross the wall at $t=t_0$.
Thus, in the case of $n=6$, for 3 of the 15 pairwise matchings the triple linking number alternates. For the remaining 12 matchings
such number is obviously preserved, so, $\sign_6(A_0)=-\sign_6(A_1)$.

In the case of $n=7$, among the seven signs of 6-subconfigurations exactly three will change: if the point that we drop is not among the four points
which become coplanar in the process of deformation $A_t$.
Thus, $\sign_7(A_0)=-\sign_7(A_1)$.
\endproof



\lemma\label{sign-alternation}
For any Cayley M-octad $\OO$ and any point $x\in\OO$ the sign $\sign(\OO,x)$ alternates after $\OO$ experiences
a 4-collision.
\endlemma

\proof
It follows from Lemma \ref{4-collision} applied to the residual 7-configuration $\OO\sm\{x\}$,
because it contains precisely one of the two complementary quadruple of points of $\OO$ involved into a 4-collision.
\endproof

\lemma\label{adjacent-signs}
If in a Cayley M-octad $\OO$ points $x_i,x_j\in\OO$ are adjacent in the graph $\G_\OO$, then $\sign(\OO,x_i)=\sign(X,x_j)$.
\endlemma

\proof
This follows from that residual 7-configurations $\OO_i=\OO\sm\{x_i\}$, $i=1,2$ are deformation equivalent through simple 7-configurations,
in which $x_1$ moves towards $x_2$ along the corresponding edge of $\G_\OO$ and the other 6 points of $\OO$ remain constant.
\endproof

\demo{Proof of Theorem \ref{chirality-of-octads} }
The signs $\sign(\OO,x)$ are the same for all $x\in\OO$ if
$\OO\in\O44$ due to Lemma \ref{adjacent-signs}, because the graph $\G_\OO$ is connected.
Since by Corollary \ref{total-adjacency} we may reach any
coarse deformation classes $\O{\a}{\b}$ through wall-crossing beginning from $\O44$,
Lemma \ref{sign-alternation} implies that for any regular M-octad $\OO$ the signs $\sign(\OO,x)$ are
the same for all eight points of $\OO$.
\qed\enddemo

By Theorem \ref{chirality-of-octads} each coarse deformation class $\O{\a}{\b}$ is the union of two pure deformation
classes distinguished by the common sign of vertices. According to this sign, we will denote the corresponding pair of deformation classes
$\O{\a}{\b}^+$ and $\O{\a}{\b}^-$.

\corollary
Any pair of deformation classes $\O{\a}{\b}^\pm$ are related by several wall-crossings.
\endcorollary

\proof
It follows from Corollary \ref{total-adjacency} and adjacency between classes $\O23^+$ and $\O23^-$ realized by
a black-edge-move in the theta-diagram $D_{23}$, for a black edge whose endpoints are of different colors.
\endproof

\section{Real monodromy groups}

\subsection{The Complex monodromy}
Consider a regular Cayley Octad $\OO$ with its Hessian quartic $\H$ and spectral theta-characteristic $\theta\in\Theta_0(\H)$.
Isomorphism $\Phi_H\:\LLL(\OO)\cong H_1(\H;\Z/2)$ due to Proposition \ref{identification} and 
relation (3) in Sect. \ref{theta-quadratic}
induces a homomorphism $\PH:S(\OO)\to \Aut(H_1(\H;\Z/2),q_\theta)$ from the permutation group $S(\OO)\cong S_8$
acting on $\OO$ to the group of automorphisms of
$H_1(\H;\Z/2)$ preserving the quadratic function $q_\theta$ associated to $\theta$
(and thus, preserving $\Z/2$-valued intersection form in $H_1(\H;\Z/2)$).

\proposition\label{representation}\!\!$($cf.\,\cite{GH}, Proposition 2.1$)$
$\PH$ is a group isomorphism.
\endproposition

\proof
Proposition \ref{identification} immediately implies that $\PH$ is monomorphic, so, it is left to
to verify that $|\Aut(H_1(\H;\Z/2),q_\theta)|=8!$.
First, we find the order of the automorphism group $\Aut(H_1(\H;\Z/2))$ preserving just the intersection form,
by counting the number of symplectic bases $e_1,f_1,e_2,f_2,e_3,f_3$ in $H_1(\H;\Z/2)$.
This number is $(63\cdot32)(15\cdot8)(3\cdot2)=36(8!)$, where, $63\cdot32$ counts the number of pairs $e_1,f_1$, etc.
Group  $\Aut(H_1(\H;\Z/2))$ acts transitively on the set $\Theta_0(\H)$ (identified with the set of quadratic function with even
Arf invariant). Since $\Aut(H_1(\H;\Z/2),q_\theta)$ is the stabilizer of $\theta$, its order is $8!$, since
$|\Theta_0(\H)|=36$.
\endproof

%
%

\subsection{The monodromy groups of regular M-octads}
By the {\it real monodromy group}, $\Aut_\R(\OO)$, of a real regular Cayley M-octad $\OO$ we mean
the subgroup of $S(\OO)$ that is the image under
the monodromy homomorphism $\pi_1(\O{\a}{\b}^\pm,\OO)\to S(\OO)$, where $\O{\a}{\b}^\pm$ is the deformation component
of M-octads containing $X$.

\theorem\label{real-monodromy}
For any regular Cayley M-octad $\OO$ the group $\Aut_\R(\OO)$ is isomorphic to a subgroup of $S_4$ formed by
the permutations of the ovals $\T_i$, $0\le i\le3$, of the associated Hessian quartic $\H$ which preserve their
colors as well as the colors of bridges $b_{ij}$, $0\le i<j\le3$ under the induced permutation of them.
\endtheorem


\proof
By Proposition \ref{oval-bridge-bases} a permutation of ovals of M-quartic $\H$ determines the corresponding permutation of bridges
and thus, an automorphism of $H_1(\H;\Z/2)$ which is necessarily real and non-trivial for a non-trivial permutation of ovals.
By Proposition \ref{representation}, an automorphism of $H_1(\H;\Z/2)$ determines a permutation of $\OO$ provided it preserves
the quadratic function $q_\theta$ associated to the spectral theta-characteristic $\theta$, which is equivalent to preserving
the coloring of ovals and edges on the theta-diagram of $(\H,\theta)$. The permutation of $\OO$ that we obtain is represented by a
real monodromy, as it follows from Theorem \ref{real-Dixon}.
\endproof

\corollary
The list of groups $\Aut_\R(\OO)$ for $\OO\in\O{\a}{\b}$, $a\in\{0,2,4\}$ and $b\in\{0,3,4\}$, $(\a,\b)\ne(2,3)$
is as presented in the Table 2 of Fig.\,1.
\endcorollary

\proof
Theorem \ref{real-monodromy} reduces finding of $\Aut_\R(\OO)$ to an trivial analysis of symmetries of
the eight theta-diagrams on Fig.\, 2.
\endproof

\section{Concluding remarks}

\subsection{Deformation classes of marked Cayley M-octads}
As an application of analysis of the real monodromy action in the last Section, we can now easily upgrade
our deformation classification of Cayley M-octad to more refined classification of M-octads with a marked point.
Namely, each coarse deformation class $\O{\a}{\b}$ gives a number of coarse deformation classes or marked M-octad which is obviously equal to
the number of orbits of the real monodromy action on $\OO\in \O{\a}{\b}$.
The numbers of orbits are indicated in the Table 3 on Figure 1 and there sum, 14, is the number of classes of marked M-octads.

Furthermore, dropping the chosen point in a marked regular M-octad gives a configuration of 7 points in $P^3_\R$.
On the other hand, the central projection from this marked point send the others into a configuration of 7 points in $P^2_\R$, 
 none 3 of which is collinear and no 6 is coconic.  Following \cite{GH} we call such 7-configurations {\it typical}.
These two 7-configurations, planar and spacial, are related by the {\it Gale duality} (see \cite{GH}, Sect.\,7),
and knowing one of them (up to projective equivalence), one can recover the whole Cayley octad
(cf., Sect.\,\ref{Aronhold-section}).

 In particular the 14 coarse deformation classes of marked regular Cayley M-octads correspond to the 14 deformation classes
of planar typical 7-configurations that were analyzed in \cite{FZ}.

\subsection{Theta-characteristics for real (M-1)-quartics}
In the case of regular (M-1)-octads the corresponding Hessian curves $\H$ are (M-1)-quartics, by Proposition \ref{discrepancy-correspondence}.
The set of ovals $\T_1,\T_2,\T_3$ of $\H$
is connected pairwise by bridges: each pair $\T_i$, $\T_j$ by two bridges denoted $b_{ij}^\pm$:
in accord with two kinds of line segments connecting a pair of points in $\Rp2$, see Fig.\,4 for the corresponding theta-diagrams.
It is not difficult to show that $[b_{ij}^+]+[b_{ij}^-]=[\H_\R]$, and so, $q_\theta$ takes opposite values on $[b_{ij}^+]$ and $[b_{ij}^-]$.
A further analysis shows that
there exist two deformation classes of theta-diagrams for even spin (M-1)-quartics, which we presented on Fig.\,4
together with the corresponding graphs $\G_X$ (which have six vertices for (M-1)-octads) and real monodromy groups
$\Aut_\R(\OO)$ (the definitions and proofs are analogous to the case of M-octads).

In particular, we see that {\it there exist precisely two coarse deformation classes of regular (M-1)-octads described by
theta-diagrams and graphs $\G_X$ on Fig.\,4}.
By Proposition \ref{6-7-chirality}, pure deformation classification of regular (M-1)-octads gives then 4 deformation classes.

\midinsert
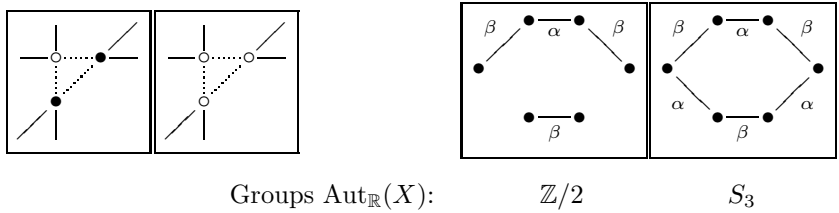
\topcaption{Fig.\,4. Theta-diagrams and graphs $\G_\OO$ for (M-1)-octads $\OO$}\endcaption
\hskip1mm
\hbox{\vbox{\hsize22.5mm
$\hbox{\boxed{\xymatrix@1@=10pt{
{}&&&{}\\
&\circ\ar@{-}[u]\ar@{-}[l]\ar@{.}[r]\ar@{.}[d]&\bullet\ar@{.}[r]\ar@{-}[ur]\ar@{-}[r]&\\
&\bullet\ar@{-}[dl]\ar@{-}[d]\ar@{.}[ur]\\
{}&&&{}
}}}$
}
\hskip-4mm
\vbox{\hsize45mm
$\boxed{\xymatrix@1@=10pt{
{}&&&{}\\
&\circ\ar@{-}[u]\ar@{-}[l]\ar@{.}[r]\ar@{.}[d]&\circ\ar@{.}[r]\ar@{-}[ur]\ar@{-}[r]&\\
&\circ\ar@{-}[dl]\ar@{-}[d]\ar@{.}[ur]\\
{}&&&{}
}}$
}\vbox{\hsize25mm
\hbox{$\boxed{\xymatrix@1@=12pt{
&\bullet\ar@{-}[r]_{\a}&\bullet\\
\bullet\ar@{-}[ur]^{\b}&&&\bullet\ar@{-}[ul]_{\b}\\
&\bullet\ar@{-}[r]_{\b}&\bullet}}$}
}
\hskip-1mm
\vbox{\hsize30mm
\hbox{$\boxed{\xymatrix@1@=12pt{
&\bullet\ar@{-}[r]_{\a}&\bullet\\
\bullet\ar@{-}[ur]^{\b}\ar@{-}[dr]_{\a}&&&\bullet\ar@{-}[ul]_{\b}\ar@{-}[dl]^{\a}\\
&\bullet\ar@{-}[r]_{\b}&\bullet}}$}
}
}
\vskip3mm
\hskip35mm {Groups $\Aut_\R(\OO)$:\hskip14mm $\Z/2$ \hskip18mm $S_3$}
\endinsert

\subsection{Odd spin quartics}
An odd theta-characteristic, $\theta$, on a non-singular quartic $C$ defines a cubic surface $Z$ with a point $z\in Z$,
so that $C$ is projectively equivalent to the critical locus of the central projection of $Z$ to $P^2$ from the center $z$,
and $\theta$ is represented by the tangent plane to $Z$ at $z$ which is projected to a bitangent to $C$.

The three deformation classes of odd spin M-quartics $(C,\theta)$ found in Proposition \ref{codes-of-theta-orbits}
corresponds to the three types of real points on a real nonsingular M-cubic surface $Z\subset P^3$.
Recall that $Z$ is {\it M-cubic} if and only if
all 27 lines on it are real, and the complement of these lines in $Z_\R$ split into polygonal regions,
which  can be triangular, quadrilateral, or pentagonal. These three types of regions correspond to the above three types of $(C,\theta)$,
according to the total number of black ovals and bridges on an odd theta-diagram, which can be respectively 3,4 or 5 (see the three rightmost diagrams of Fig.\,3).

\enddocument